\documentclass{ifacconf} 

\usepackage{natbib}

\counterwithin*{section}{part}

\usepackage{amsmath,amssymb}
\usepackage{graphicx}

\newenvironment{proof}{\begin{pf}}{\hspace{\fill} \qed \end{pf}}

\usepackage[noabbrev,capitalize]{cleveref}

\newcommand{\bbR}{\mathbb{R}}

\newcommand{\bbX}{\mathbb{X}}
\newcommand{\bbU}{\mathbb{U}}
\newcommand{\bbZ}{\mathbb{Z}}
\newcommand{\bbW}{\mathbb{W}}

\newcommand{\bbB}{\mathbb{B}}

\newcommand{\bfu}{\mathbf{u}}
\newcommand{\bfw}{\mathbf{w}}

\newcommand{\calU}{\mathcal{U}}
\newcommand{\calX}{\mathcal{X}}

\newcommand{\calK}{\mathcal{K}}

\newcommand{\norm}[1]{\lVert #1 \rVert}

\crefname{equation}{}{}
\crefname{Equation}{}{}
\crefname{assum}{Assumption}{Assumptions}
\crefname{thm}{Theorem}{Theorems}
\crefname{defn}{Definition}{Definitions}

\begin{document}
	
\begin{frontmatter}

\title{Inherently Robust Economic Model Predictive Control Without 
	Dissipativity}

\author{Robert D. McAllister}
\address{Delft Center for Systems and Control, Delft University of 
Technology, Delft, NL (e-mail: r.d.mcallister@tudelft.nl)}

\begin{abstract}
    We establish sufficient conditions for the terminal cost and constraint 
    such that economic model predictive control (MPC) is robustly recursively 
    feasible and economically robust to small disturbances without any 
    assumptions of dissipativity. Moreover, we demonstrate that these 
    sufficient conditions can be satisfied with standard design methods. A 
    small example is presented to illustrate the inherent robustness of 
    economic MPC to small disturbances. 
\end{abstract}
\end{frontmatter}
\section{Introduction}

For successful implementation, model predictive control (MPC) must be robust to 
disturbances such that arbitrarily small perturbations and modeling errors 
produce similarly small deviations in performance. For setpoint tracking 
problems, in which the stage cost is positive definite with respect to this 
setpoint, nominal MPC ensures a nonzero margin of inherent robustness to 
disturbances and prediction errors 
\citep{yu:reble:chen:allgower:2014,allan:bates:risbeck:rawlings2017}. 

For economic MPC problems, in which the stage cost is not necessarily positive 
definite with respect to a setpoint, asymptotic stability of a specific 
steady-state can be guaranteed via an assumption of strict dissipativity 
\citep{diehl:amrit:rawlings:2010,angeli:amrit:rawlings:2011,amrit:rawlings:angeli:2011}.
Without terminal costs/constraints, strict dissipativity is used to 
establish practical asymptotic stability of an optimal, but potentially 
unknown, steady-state 
\citep{grune:2013,faulwasser:bonvin:2015,grune:muller:2016}. 
Gradient-correcting end penalties can also be used to ensure asymptotic 
stability without terminal constraints for sufficiently long horizons 
\citep{zanon:faulwasser:2018}. These results all rely on 
strict dissipativity to ``rotate'' the economic MPC problem such that tools and 
results from tracking MPC can be applied. Thus, extending the 
guarantees of inherent robustness for tracking MPC to strictly dissipative 
economic MPC formulations is straightforward. Verifying strict dissipativity of 
a specific steady state is nontrivial, but there are systematic methods 
available that use sum-of-squares techniques 
\citep{pirkelmann:angeli:grune:2019}. 

Unfortunately, this assumption of strict dissipativity does not always hold for 
a specific steady state. Thus, purely economic stage costs are often 
modified to ensure strict dissipativity and thereby compromise potential 
economic gains for guaranteed stability of a chosen steady-state target. For 
example, \citet{zanon:gros:diehl:2016} present an approach to design an 
asymptotically stable tracking MPC formulation that is locally equivalent to 
the original economic MPC problem. To encourage steady-state operation, 
\citet{alamir:pannocchia:2021} include a rate-of-change penalty for the 
open-loop state trajectory.

In some applications, however, economic performance is more 
important than tracking a specific setpoint for the system and dynamic 
operation of the system may improve economic performance (e.g., production 
scheduling, HVAC, energy systems). While an optimal (periodic) 
trajectory can be used instead of a steady-state, verifying strict 
dissipativity for this trajectory is quite difficult 
\citep{muller:grune:2016,grune:pirkelmann:2020}. For polynomial optimal control 
problems, sum-of-squares techniques can also be used to systematically check 
for strictly dissipative periodic trajectories 
\citep{berberich:kohler:allgower:muller:2020}. However, these methods can 
be computationally expensive and a strictly dissipative periodic trajectory may 
not exist for some optimal control problems of interest. 

Without dissipativity, little is known about the inherent 
robustness of economic MPC and the results from tracking MPC are not directly 
applicable. In previous work, we addressed the inherent robustness of economic 
MPC subject to large and infrequent disturbances, but avoided the issue of 
recursive feasibility by assuming that the economic 
MPC problem was recursively feasible by design 
\citep{mcallister:rawlings:2023}. Robust EMPC formulations guarantee 
recursive feasibility via constraint tightening and can thereby establish 
robust performance guarantees \citep{bayer:lorenzen:muller:allgower:2016, 
dong:angeli:2020,schwenkel:kohler:muller:allgower:2020}. However, these 
constraint-tightening procedures are nontrivial for nonlinear systems. 
To the best of our knowledge, there are no results that guarantee the inherent 
robustness of economic MPC without either constraint tightening or 
dissipativity. 

\textbf{Contribution:} In this work, we focus on economic MPC problems in which 
a reasonable steady 
state for the system is available to serve as a baseline, but operating near 
this steady-state is not required. The main contribution of this work is a set 
of requirements for the terminal cost and constraint (\cref{as:term}) that are 
sufficient to guarantee that economic MPC is robustly recursively feasible and 
inherently robust to small disturbances in terms of economic performance 
(\cref{thm:robust}) without any assumptions of dissipativity or constraint 
tightening. We then demonstrate how standard procedures 
to construct terminal costs and constraints for economic MPC, first presented 
in \citet{amrit:rawlings:angeli:2011}, can satisfy these requirements. We 
conclude with a small example to demonstrate the implications of this analysis.

\textbf{Notation:} Let $\bbR$ denote the reals with 
subscripts and superscripts denoting the restrictions and dimensions (e.g., 
$\bbR_{\geq 0}^n$ for nonnegative reals of dimension $n$). Let $|\cdot|$ denote 
Euclidean norm. 
Let $\varepsilon\bbB:=\{x\in\bbR^n\mid |x|\leq \varepsilon\}$. The function 
$\alpha:\bbR_{\geq 
0}\rightarrow\bbR_{\geq 0}$ is in class $\calK$ if 
it is continuous, strictly increasing, and $\alpha(0)=0$. The function 
$\alpha(\cdot)$ is in class $\calK_{\infty}$ if  
$\alpha(\cdot)\in\calK$ and $\lim_{s\rightarrow\infty}\alpha(s)=\infty$. Let 
$\nabla f(x)$ ($\nabla^2f(x)$) denote the gradient (Hessian) of $f(\cdot)$ at 
$x$. 

\section{Problem formulation and preliminaries}
We consider the discrete-time system
\begin{equation*}
    x^+ = f(x,u,w) \qquad f:\bbX\times\bbU\times\bbW \rightarrow \bbX
\end{equation*}
in which $x\in\bbX\subseteq\bbR^n$ is the state, $u\in\bbU\subseteq\bbR^m$ is 
the input, $w\in\bbW\subseteq\bbR^{q}$ is the disturbance, and $x^+\in\bbX$ is 
the successor state. In the economic MPC problem, we use the nominal model:
\begin{equation}\label{eq:fnom}
    x^+ = f(x,u,0)
\end{equation}
For the horizon $N\geq 1$, let $\hat{\phi}(k;x,\bfu)$ denote the state of the 
dynamical system in \cref{eq:fnom} at time $k\in\{0,\dots,N\}$ given the 
initial state $x\in\bbX$ and input trajectory 
$\bfu:=(u(0),\dots,u(N-1))\in\bbU^N$. 


Fundamental physical limits of the system (e.g., temperature cannot be below 0 
K) can be enforced via the function $f(\cdot)$ and thereby represented by 
$\bbX$, i.e., the range of $f(\cdot)$. In many applications of economic MPC, 
\emph{desired} state (mixed) constraints are also relevant, i.e., we want 
\begin{equation}\label{eq:hard_con}
    (x,u) \in \bbZ_g := \{ (x,u)\in\bbX\times\bbU \mid g(x,u) \leq 0\}
\end{equation}
for some continuous function $g:\bbX\times\bbU\rightarrow\bbR$. In general, 
nominal MPC cannot guarantee that this general class of constraints is 
satisfied for a perturbed system. Thus, enforcing $(x,u) \in\bbZ_g$ as a 
\emph{hard} constraint in the MPC optimization problem can easily lead to 
infeasible optimization problems and therefore undefined control laws. Instead, 
these constraints are \textit{softened} via a penalty function and the stage 
cost becomes 
\begin{equation}\label{eq:ell}
    \ell(x,u):=\ell_e(x,u) + \lambda \max\{g(x,u),0\}
\end{equation}
in which $\ell_e:\bbX\times\bbU\rightarrow\bbR$ is the economic cost and 
$\lambda >0$ defines the weighting of the penalty function.\footnote{For 
numerical optimization, $\max\{g(x,u),0\}$ can be rewritten via a slack 
variable $s\geq 0$ with the constraint $s\geq g(x,u)$.} 


While we do not permit a general class of hard state constraints, we do enforce 
a hard constraint on the terminal state in the prediction horizon, i.e., 
\begin{equation*}
    \hat{\phi}(N;x,\bfu)\in\bbX_f
\end{equation*}
in which $\bbX_f\subseteq\bbX$. Unlike the desired constraint $\bbZ_g$, 
however, this terminal constraint must satisfy specific assumptions (e.g., 
\cref{as:term}) and is only enforced on the \emph{final} state in the open-loop 
trajectory. We also define the terminal cost $V_f:\bbX\rightarrow\bbR$. 

With only input constrains and a terminal constraint, we denote the set of admissible control trajectories as
\begin{equation*}
    \calU(x) := \left\{\bfu\in\bbU^N \ \middle| \ \hat{\phi}(N;x,\bfu)\in\bbX_f \right\}
\end{equation*}
and the set of all feasible initial states as
\begin{equation*}
    \calX := \{x\in\bbX \mid \calU_N(x)\neq\emptyset\}
\end{equation*}
We define cost function
\begin{equation*}
    V(x,\bfu) := \sum_{k=0}^{N-1}\ell(\phi(k;x,\bfu),u(k)) + V_f(\phi(N;x,\bfu))
\end{equation*}
The optimal MPC problem is then
\begin{equation*}
    \ V^0(x) := \min_{\bfu\in\calU(x)}V(x,\bfu) 
\end{equation*}
and the optimal solution(s) are defined as 
$\bfu^0(x):=\arg\min_{\bfu\in\calU(x)}V(x,\bfu)$. 

The control law $\kappa:\calX\rightarrow\bbU$ is defined as 
$\kappa(x):=u^0(0;x)$ in which $u^0(0;x)$ is the first input in the trajectory 
$\bfu^0(x)$. The closed-loop systems is therefore
\begin{equation}
    x^+ = f(x,\kappa(x),w) \label{eq:fcl}
\end{equation}
Let $\phi(k;x,\bfw_{\infty})$ denote the closed-loop state of \cref{eq:fcl} at 
time $k\geq 0$ given the initial state $x\in\calX$ and disturbance trajectory 
$\bfw_{\infty}=(w(0),w(1),\dots)\in\bbW^{\infty}$. Let 
$\norm{\bfw_{\infty}}:=\sup_{k\geq 
0}|w(k)|$ and define the following terms.
\begin{defn}[Positive invariant]
    A set $X$ is positive invariant for $x^+= f(x)$ if $f(x)\in X$ for all 
    $x\in X$
\end{defn}

\begin{defn}[Robustly positive invariant]
    A set $X$ is robustly positive invariant (RPI) for $x^+= f(x,w)$, $w\in W$ 
    if $f(x,w)\in X$ for all $x\in X$ and $w\in W$. 
\end{defn}

\begin{defn}[Lyapunov function]
    The function $V:X\rightarrow\bbR_{\geq 0}$ is a Lyapunov function for 
    $x^+=f(x)$ on the positive invariant set $X\subseteq\bbR^n$ with respect to 
    the steady-state $x_s$ if there exist 
    $\alpha_1(\cdot),\alpha_2(\cdot),\alpha_3(\cdot)\in\calK_{\infty}$ such that
    \begin{subequations}
    \begin{gather}
        \alpha_1(|x-x_s|) \leq V(x) \leq \alpha_2(|x-x_s|) \\
        V(f(x)) \leq V(x) - \alpha_3(|x-x_s|)\label{eq:lyap_dec_def}
    \end{gather}
    \end{subequations}
    for all $x\in X$. 
\end{defn}

\section{Main technical results}

To construct a terminal cost and constraint, we first select a high-quality 
(low cost) steady-state pair $(x_s,u_s)\in\bbZ_g$ for the system to serve as a 
reference. For example, 
\begin{equation*}
    (x_s,u_s) \in \arg\min\left\{\ell(x,u) \ \middle| \ x = f(x,u,0), \ (x,u)\in\bbZ_g\right\}
\end{equation*}
We consider the following standard regularity assumption. 

\begin{assum}[Continuity and closed-sets]\label[assum]{as:cont}
    The system $f:\bbX\times\bbU\times\bbW\rightarrow\bbX$ and stage cost $\ell:\bbX\times\bbU\rightarrow\bbR$ are continuous. The set $\bbU$ is compact, $\bbX$ is closed, and $\calX$ is bounded. The pair $(x_s,u_s)\in\bbX\times\bbU$ satisfies $x_s=f(x_s,u_s,0)$.  
\end{assum}

\begin{rem}[Bounded $\calX$]
    If $\bbX_f$ and $\bbU$ are compact and the set $f^{-1}(S):=\{(x,u)\in\bbX\times\bbU \mid f(x,u,0)\in S\}$ is bounded for all bounded $S\subseteq\bbX$, then $\calX$ is also bounded \cite[Prop. 2.10d]{rawlings:mayne:diehl:2020}. The requirement that  $f^{-1}(S)$ is bounded for all bounded $S\subseteq\bbX$ is a mild requirement if $f(\cdot)$ is the discrete time version of a continuous system (see \cite[p. 111]{rawlings:mayne:diehl:2020}). 
\end{rem} 

\subsection{Nominal performance guarantees}

The standard assumption for the terminal cost and constraint is as follows. 
\begin{assum}[Standard terminal cost and 
constraint]\label[assum]{as:term_standard}\ \\
	The terminal set $\bbX_f\subseteq\bbX$ is compact and $x_s\in\bbX_f$. The 
	terminal cost $V_f:\bbX\rightarrow\bbR$ is continuous. There exists a 
	terminal control law $\kappa_f:\bbX_f\rightarrow\bbU$ such that 
	$f(x,\kappa_f(x),0)\in\bbX_f$ and
    \begin{equation*}
        V_f(f(x,\kappa_f(x),0)) \leq V_f(x) -\ell(x,\kappa_f(x)) + \ell(x_s,u_s)
    \end{equation*}
    for all $x\in\bbX_f$.
\end{assum}

We note that selecting a terminal cost and constraint according to \cref{as:term_standard} is simple. For example, one can choose $\bbX_f=\{x_s\}$, $V_f(x_s)=0$, and $\kappa_f(x_s)=u_s$ to satisfy \cref{as:term_standard}. This assumption is sufficient to establish the following nominal performance guarantee.

\begin{thm}[Nominal performance 
\citep{angeli:amrit:rawlings:2011}]\label[thm]{thm:nom}
    Let \cref{as:cont,as:term_standard} hold. Then we have that $\calX$ is positive invariant for \cref{eq:fcl} with $w=0$ and 
    \begin{equation*}
        \limsup_{T\rightarrow\infty}\frac{1}{T}\sum_{k=0}^{T-1} \ell(x(k),u(k)) \leq \ell(x_s,u_s)
    \end{equation*}
    in which $x(k)=\phi(k;x,\mathbf{0})$, $u(k)=\kappa(x(k))$ for all 
    $x\in\calX$. 
\end{thm}

\cref{thm:nom} ensures that the long-term performance, in terms of the stage 
cost, for the closed-loop system generated by economic MPC is no worse than the 
steady-state reference $\ell(x_s,u_s)$. Note that \cref{thm:nom} does not 
guarantee asymptotic stability of this steady-state. The closed-loop 
trajectory may instead follow a periodic (or an aperiodic) orbit. 
Nonetheless, the guarantee in \cref{thm:nom} ensures that economic MPC does not 
generate unnecessarily poor closed-loop performance, even for small values of 
$N$. This guarantee is particularly important when the dynamics $f(\cdot)$ are 
complicated and/or high dimensional, in which case using long horizons in the 
optimization problem may not be tractable. 

\subsection{Inherent robustness}

Unfortunately, \Cref{as:term_standard} may not provide any inherent robustness 
to the controller (see \citet{grimm:messina:tuna:teel:2004} for examples). In 
particular, 
arbitrarily small disturbances $w$ may render the economic MPC problem 
infeasible if a terminal \emph{equality} constraint is used, i.e., 
$\bbX_f=\{0\}$. To ensure a nonzero margin of inherent robustness for the 
economic MPC controller, we use a stronger version of \cref{as:term_standard}. 
Specifically, we assume that the terminal control law $\kappa_f(\cdot)$ is 
\emph{stabilizing} and the terminal set $\bbX_f$ is defined as a sublevel set 
of a local Lyapunov function.  

\begin{assum}[Robust terminal cost and constraint]\label[assum]{as:term} \ \\
    There exists a terminal control law $\kappa_f:\bbX_s\rightarrow\bbU$ and 
    continuous Lyapunov function $V_s:\bbX\rightarrow\bbR_{\geq 0}$ for the 
    system $x^+=f(x,\kappa_f(x),0)$ in the positive invariant set 
    $\bbX_{s}\subseteq\bbX$ with respect to the steady-state 
    $x_s\in\bbX_{s}$. The terminal set is defined as 
    \begin{equation}\label{eq:termcon}
        \bbX_f := \{x\in\bbX \mid V_{s}(x) \leq \tau\}
    \end{equation}
    with $\tau >0$ chosen such that $\bbX_f\subseteq\bbX_s$. The function $V_f:\bbX\rightarrow\bbR$ is continuous and satisfies 
    \begin{equation}\label{eq:termcostdec}
        V_f(f(x,\kappa_f(x),0)) \leq V_f(x) - \ell(x,\kappa_f(x)) + \ell(x_s,u_s)
    \end{equation}
    for all $x\in\bbX_f$. 
\end{assum}

Note that the terminal constraint in \cref{eq:termcon} is similar to the 
terminal constraint required to establish the robustness of nominal MPC with a 
positive definite stage cost (compare with 
\cite{yu:reble:chen:allgower:2014,allan:bates:risbeck:rawlings2017}). Unlike 
nominal MPC, however, we do not set the terminal cost equal to this local 
Lyapunov function, i.e., we allow for $V_s(x)\neq V_f(x)$. The terminal cost is 
instead defined to satisfy the cost decrease condition in 
\cref{eq:termcostdec}. Thus, the key novel feature of \cref{as:term} is that 
this assumption partially separates the design of the terminal cost and 
constraint. 

\begin{rem}[Reduction to tracking MPC]
    If there exists $\alpha_{\ell}(\cdot)\in\calK$ such that $\ell(x,u) \geq \alpha_{\ell}(|x-x_s|)$ and $\ell(x_s,u_s)=0$, we recover a tracking MPC formulation. In this case, \cref{eq:termcostdec} is equivalent to \cref{eq:lyap_dec_def} and $V_f(\cdot)$ is effectively required to be a local Lyapunov function for $x^+=f(x,\kappa_f(x),0)$. Thus, choosing $V_s(x)=V_f(x)$ satisfies \cref{as:term}. 
\end{rem}

We now establish the main result of this paper.
\begin{thm}[Inherently robust economic performance]\label[thm]{thm:robust}
    \ \\ Let \cref{as:cont,as:term} hold. Then there exist $\delta>0$ and 
    $\gamma(\cdot)\in\calK$ such that $\calX$ is RPI 
    for \cref{eq:fcl} with $w\in\{w\in\bbW\mid |w|\leq\delta\}$ and 
    \begin{equation}\label{eq:performance_robust}
        \limsup_{T\rightarrow\infty} \frac{1}{T}\sum_{k=0}^{T-1} \ell(x(k),u(k)) \leq \ell(x_s,u_s) + \gamma\left(\norm{\bfw_{\infty}}\right)
    \end{equation}
    in which $x(k)=\phi(k;x,\bfw_{\infty})$ and $u(k)=\kappa(x(k))$ for all $x\in\calX$ and $\bfw_{\infty}\in\{w\in\bbW\mid |w|\leq\delta\}^{\infty}$. 
\end{thm}
\noindent \Cref{thm:robust} ensures that arbitrarily small disturbances:
\begin{enumerate}
    \item Do not render the EMPC optimization problem infeasible ($\calX$ is RPI)
    \item Produce a similarly small degradation in the nominal performance guarantee, given by $\gamma(\lVert \bfw_{\infty}\rVert)$
\end{enumerate} 
We note that the function $\gamma(\cdot)$ is typically too conservative to 
provide useful quantitative information about the closed-loop system. 
Nonetheless, \cref{thm:robust} ensures that such a function exists and thereby 
prevents arbitrarily poor closed-loop performance systems with for small 
disturbances. In other words, the EMPC controller is not fragile in a practical 
setting.  

\begin{pf}[Proof of \cref{thm:robust}]
    Since $f(\cdot)$ is continuous and $\calX\times\bbU$ is bounded, we have from \cite[Prop. 20]{allan:bates:risbeck:rawlings2017} that there exists $\sigma_f(\cdot)\in\calK_{\infty}$ such that
    \begin{equation*}
        |f(x,u,w) - f(x,u,0)| \leq \sigma_f(|w|)
    \end{equation*}
    for all $(x,u)\in\calX\times\bbU$ and $w\in\bbW$. 
    Since $f(\cdot)$ and $V_s(\cdot)$ are continuous, 
    $V_s(\hat{\phi}(N;x,\bfu))$ is continuous. From \cite[Prop. 
    20]{allan:bates:risbeck:rawlings2017}, there exists 
    $\sigma_s(\cdot)\in\calK_{\infty}$ such that
    \begin{equation*}
        \rvert V_s(\hat{\phi}(N;x_1,\bfu)) - V_s(\hat{\phi}(N;x_2,\bfu))\lvert \leq \sigma_s\left(\lvert x_1-x_2\rvert\right))
    \end{equation*}
    for all $x_1\in\bbX$, $x_2\in\{f(x,u,0) \mid x\in\calX, \ u\in\bbU\}$, and 
    $\bfu\in\bbU^N$. Substituting $x_1=f(x,\kappa(x),w)$ and 
    $x_2=f(x,\kappa(x),0)$, we have
    \begin{multline}
        \rvert V_s(\hat{\phi}(N;f(x,\kappa(x),w),\bfu)) - V_s(\hat{\phi}(N;f(x,\kappa(x),0),\bfu))\lvert \\ \leq \sigma_s\left(\lvert f(x,\kappa(x),w)-f(x,\kappa(x),0)\rvert\right) \leq \sigma(|w|) \label{eq:Vs_sigma}
    \end{multline}
    in which $\sigma(\cdot) := \sigma_s\circ\sigma_f(\cdot)\in\calK_{\infty}$ for all $x\in\calX$, $\bfu\in\bbU^N$ and $w\in\bbW$. 

    For any $x\in\calX$ and $w\in\bbW$, define $\bfu^0=\bfu^0(x)$, 
    $x^+=f(x,\kappa(x),w)$, $\hat{x}^+=f(x,\kappa(x),0)$ and $x(N) := 
    \hat{\phi}(N;x,\bfu^0(x))$. With the terminal control law in 
    \cref{as:term}, we construct the candidate trajectory
    \begin{equation*}
        \tilde{\bfu}^+ := (u^0(1),\dots,u^0(N-1),\kappa_f(x(N)))
    \end{equation*}
    Denote $\hat{x}^+(N):=\hat{\phi}(N;\hat{x}^+,\tilde{\bfu}^+)$, 
    $x^+(N):=\hat{\phi}(N;x^+,\tilde{\bfu}^+)$. Note that $\hat{x}^+(N) = 
    f(x(N),\kappa_f(x(N)),0)$ by the definition of $\tilde{\bfu}^+$. 
    From \cref{as:term} and \cref{eq:lyap_dec_def}, there exists $\alpha_3(\cdot)\in\calK_{\infty}$ such that
    \begin{equation*}
        V_s(\hat{x}^+(N)) \leq  V_s(x(N)) -\alpha_3(|x(N)-x_s|)
    \end{equation*}
    and by using \cref{eq:Vs_sigma} we have 
    \begin{equation*}
        V_s(x^+(N)) \leq V_s(x(N))  - \alpha_3(|x(N)-x_s|) + \sigma(|w|)
    \end{equation*}
    Moreover, there exist $\alpha_1(\cdot),\alpha_2(\cdot)\in\calK_{\infty}$ such that
    \begin{equation*}
        \alpha_1(|x(N)-x_s|)\leq V_s(x(N)) \leq \alpha_2(|x(N)-x_s|)
    \end{equation*}
    Since $x(N)\in\bbX_f$, we have that $V_s(x(N))\leq \tau$. If $V_s(x(N)) 
    \leq \tau/2$ and $|w|\leq \sigma^{-1}(\tau/2)$, then
    \begin{align*}
        V_s(x^+(N)) & \leq \tau/2  - \alpha_3(|x(N)-x_s|) + \sigma(|w|) \\
        & \leq \tau/2 + \sigma(|w|) \leq \tau
    \end{align*}
    and therefore $x^+(N)\in\bbX_f$. If $V_s(x(N))\in[\tau/2,\tau]$ and $|w|\leq \sigma^{-1}(\alpha_3(\alpha_2^{-1}(\tau/2)))$, then
    \begin{align*}
        V_s(x^+(N)) & \leq \tau  - \alpha_3(|x(N)|) + \sigma(|w|) \\
        & \leq \tau - \alpha_3(\alpha_2^{-1}(V_s(x(N)))) + \sigma(|w|) \\
        & \leq \tau - \alpha_3(\alpha_2^{-1}(\tau/2)) + \sigma(|w|) \leq \tau
    \end{align*}
    and therefore $x^+(N)\in\bbX_f$. Thus, we define 
    $\delta:=\min\{\sigma^{-1}(\tau/2), 
    \sigma^{-1}(\alpha_3(\alpha_2^{-1}(\tau/2)))\} > 0$.
     For all $|w|\leq \delta$, we have that $x^+(N)\in\bbX_f$. Therefore, 
    $\tilde{\bfu}^+\in\calU(x^+)$ and $x^+\in\calX$, i.e., $\calX$ is RPI for 
    \cref{eq:fcl} with $|w|\leq\delta$.

    We now consider the evolution of the cost function $V(\cdot)$. From \cref{as:term}, we have that
    \begin{equation}\label{eq:V_nomcostdec}
        V(\hat{x}^+,\tilde{\bfu}^+) \leq V^0(x) - \ell(x,\kappa(x)) +\ell(x_s,u_s)
    \end{equation}
    for all $x\in\calX$. 
    Note that $V(\cdot)$ is continuous and $\hat{x}^+\in\calX$, which is a 
    compact set. From \cite[Prop. 20]{allan:bates:risbeck:rawlings2017}, there 
    exists $\sigma_{V}(\cdot)\in\calK$ such that
    \begin{equation}\label{eq:V_diff}
        \lvert V(x^+,\tilde{\bfu}^+) - V(\hat{x}^+,\tilde{\bfu}^+) \rvert \leq 
        \sigma_{V}(|x^+-\hat{x}^+|) \leq \gamma(|w|)
    \end{equation}
    in which $\gamma(\cdot) := \sigma_V\circ\sigma_f(\cdot)\in\calK$ for all $x\in\calX$ and $w\in\bbW$. By combining \cref{eq:V_nomcostdec}, \cref{eq:V_diff}, and $V^0(x^+)\leq V(x^+,\tilde{\bfu}^+)$ we have
    \begin{equation}\label{eq:V_costdec}
        V^0(x^+) \leq V^0(x) - \ell(x,\kappa(x)) +\ell(x_s,u_s) + \gamma(|w|) 
    \end{equation}
    for all $x\in\calX$ and $w\in\bbW$. 

    Choose any $x\in\calX$, $\bfw_{\infty}\in\{w\in\bbW\mid |w|\leq 
    \delta\}^\infty$. Denote $x(k)=\phi(k,x,\bfw_{\infty})$, and 
    $u(k)=\kappa(x(k))$. Since $\calX$ is RPI for the system 
    $x^+=f(x,\kappa(x),w)$ and all $|w|\leq\delta$, we have that $x(k)\in\calX$ 
    for all $k\geq 0$. Therefore, $u(k)$ is well defined. From 
    \cref{eq:V_costdec} and $|w(k)|\leq\lVert \bfw_{\infty}\rVert$, we have
    \begin{multline*}
        \ell(x(k),u(k)) \leq V^0(x(k)) - V^0(x(k+1)) \\ +\ell(x_s,u_s) + 
        \gamma(\lVert \bfw_{\infty}\rVert)
    \end{multline*}
    We sum both sides of this inequality from $k=0$ to $T\geq 1$, divide by $T$, and rearange to give
    \begin{multline*}
        \frac{1}{T}\sum_{k=0}^{T-1} \ell(x(k),u(k)) \leq \frac{V^0(x(0)) - 
        V^0(x(T))}{T} \\ + \ell(x_s,u_s) + \gamma(\lVert \bfw_{\infty}\rVert)
    \end{multline*}
    Since $\calX\times\bbU^N$ is bounded and $V(\cdot)$ is continuous, $V^0(x(0)) - V^0(x(k))$ is bounded. Thus, we take the limit supremum as $T\rightarrow\infty$ to give \cref{eq:performance_robust}.
\end{pf}

\section{Designing terminal costs and constraints}

We now demonstrate that \cref{as:term} can be satisfied by designing the 
terminal cost and constraint via methods presented in \citet[section 
4]{amrit:rawlings:angeli:2011}. We restate this approach with some 
modifications here to demonstrate consistency with \cref{as:term}. 
We subsequently assume, without loss of 
generality, that 
\begin{equation*}
    (x_s,u_s)=(0,0)
\end{equation*} 
and consider the following assumption. 
 
\begin{assum}[Stabilizable]\label[assum]{as:stab}
    The functions $f(\cdot)$ and $\ell(\cdot)$ are twice continuously 
    differentiable in the interior of $\bbZ_g$, and the linearized system 
    $x^+=Ax+Bu$ with $A:=\frac{\partial f}{\partial x}(0,0,0)$ and 
    $B:=\frac{\partial f}{\partial u}(0,0,0)$ is stabilizable. The 
    sets $\bbU$ and $\bbZ_g$ contain the origin in their interior.  
\end{assum}

We now construct the terminal constraint and cost. Choose $K\in\bbR^{m\times 
n}$ such that $A_K:=(A+BK)$ is Schur stable. For some 
$\tilde{Q}\succ 0$, define $\tilde{P}\succ 0$ to solve the Lyapunov 
equation:
\begin{equation}\label{eq:dlyap}
	A_K'\tilde{P}A_K - \tilde{P} + \tilde{Q} = 0
\end{equation}
We define the candidate Lyapunov function
\begin{equation}\label{eq:V_s}
	V_s(x):=x'\tilde{P}x
\end{equation}
and for some $\tau>0$ we define the terminal constraint
\begin{equation}\label{eq:bbX_f}
	\bbX_f := \{x\in\bbX \mid V_s(x)\leq\tau\}
\end{equation}
such that $Kx\in\bbU$ and $(x,Kx)\in\bbZ_g$ for all $x\in\bbX_f$.  

For the terminal cost, we define
\begin{equation*}
	\bar{\ell}(x) := \ell(x,Kx) - \ell(0,0)
\end{equation*}
and $q=\nabla \bar{\ell}(0)$. From \citet[Lemma 
22]{amrit:rawlings:angeli:2011}, 
there exists a symmetric (possibly indefinite) matrix 
$Q\in\bbR^{n\times n}$ such that
\begin{equation*}
	x'(Q-\nabla^2 \bar{\ell}(x))x \geq 0 \quad \forall \ x\in\bbX_f
\end{equation*}
We define the function 
\begin{equation*}
	V_e(x) := x'Px + p'x
\end{equation*}
in which $P$ and $p$ satisfy
\begin{equation}\label{eq:dlyap_cost}
	A_K'PA_K - P +  Q = 0
\end{equation}
\begin{equation}\label{eq:p_cost}
	p' = q'(I-A_K)^{-1}
\end{equation}
For $\mu\geq 0$, the terminal cost is defined as
\begin{equation}\label{eq:V_f}
	V_f(x) := \mu V_s(x) + V_e(x)
\end{equation}

\begin{rem}[Comparison to \citet{amrit:rawlings:angeli:2011}]
	In contrast to \citet[Assumption 
	19]{amrit:rawlings:angeli:2011}, we require $\ell(\cdot)$ 
	to be twice continuously differentiable only in the interior of $\bbZ_g$ to 
	permit stage costs such as \cref{eq:ell} that include softened constraints. 
	Also, unlike \citet[See below equation (22)]{amrit:rawlings:angeli:2011}, 
	we do not modify $Q$ to ensure that $P$ is positive definite. Thus, 
	$V_f(x)$ is not necessarily convex or positive definite with respect to 
	$x=-P^{-1}p$.
\end{rem}

For sufficiently small $\tau>0$ and sufficiently large $\mu\geq 0$, this 
terminal cost and constraint satisfy \cref{as:term}. 

\begin{lem}[Constructing terminal costs and constraints]\label{lem:construct}
	If \cref{as:cont,as:stab} hold, then there exist $\tau>0$ and $\mu\geq 
	0$ such that $\kappa_f(x)=Kx$, \cref{eq:V_s}, \cref{eq:bbX_f}, and 
	\cref{eq:V_f} satisfy \cref{as:term}. 
\end{lem}

\begin{proof}
	Choose $\tilde{\tau}_1>0$ such that 
	$\tilde{\bbX}_1:=\{x \in \bbX 
	\mid V_s(x) \leq \tilde{\tau}_1\}\subset\{x\mid Kx\in\bbU, \ 
	(x,Kx)\in\bbZ_g\}$. We now show that there exists some sufficiently small 
	$\tau\in(0,\tilde{\tau}_1]$ 
	such that $V_s(\cdot)$ is a Lyapunov function for the system 
	$x^+=f(x,\kappa_f(x),0)$ on the set $\bbX_f$ defined in \cref{eq:bbX_f}.
	From \cref{eq:dlyap}, there exists $c_2>0$ such that
	\begin{equation*}
		V_{s}(x^+) - V_s(x) \leq -c_2|x|^2 +  V_{s}(x^+) - V_{s}(A_Kx)
	\end{equation*}
	Define $e(x) := f(x,Kx,0) - A_Kx$ so that
	\begin{equation}\label{eq:Vs_diff_e}
		V_{s}(x^+) - V_{s}(A_Kx) = 2(A_Kx)'\tilde{P}e(x) + e(x)'\tilde{P}e(x)
	\end{equation}
	Since $f(x,Kx,0)$ is twice continuously differentiable for all 
	$x\in\tilde{\bbX}_1$ and $\tilde{\bbX}_1$ is bounded, there exists 
	$c_{\delta}$ such that $|e(x)|\leq c_{\delta}|x|^2$ for all 
	$x\in\tilde{\bbX}_1$ \citep[pg. 
	141]{rawlings:mayne:diehl:2020}. From 
	\cref{eq:Vs_diff_e} and this bound on $|e(x)|$, there exist $c_3,c_4>0$ 
	such that
	\begin{equation*}
		V_{s}(x^+) - V_{s}(A_Kx) \leq c_3|x|^3 + c_4|x|^4
	\end{equation*}
	for all $x\in\tilde{\bbX}_1$. 
	Choose $\tilde{\varepsilon}>0$ such that 
	\begin{equation*}
		c_3|x|^3 + c_4|x|^4 \leq (c_2/2)|x|^2 \quad \forall 
		x\in\tilde{\varepsilon}\bbB
	\end{equation*}
	and therefore 
	\begin{equation}\label{eq:lyap_nonlin}
		V_s(x^+) \leq V_s(x) - (c_2/2)|x|^2 
	\end{equation}
	for all $x\in\tilde{\bbX}_1\cap\tilde{\varepsilon}\bbB$. Since 
	$\tilde{P}\succ 0$, we can choose $\tau \in (0,\tilde{\tau}_1]$ such that 
	$\bbX_f=\{x\in\bbX\mid V_s(x)\leq \tau\} \subseteq 
	\tilde{\varepsilon}\bbB$. Since $\tau\leq\tilde{\tau}_1$, we also have 
	$\bbX_f\subseteq\tilde{\bbX}_1$. Thus, \cref{eq:lyap_nonlin} holds for all 
	$x\in\bbX_f$ ensuring that $\bbX_f$ is positive invariant ($V_s(x^+)\leq 
	V_s(x)$). Moreover, $V_s(x)$ is a Lyapunov function in $\bbX_s=\bbX_f$ with 
	respect to the origin $x_s=0$.

	We now address the terminal cost. From \citet[Lemma 
	23]{amrit:rawlings:angeli:2011} and the 
	definition of $V_s(x)$ and $V_e(x)$, we have
	\begin{equation*}
		V_f(x^+) - V_f(x) \leq -\bar{\ell}(x) - \frac{\mu}{2}c_2|x|^2 + 
		V_{e}(x^+) - V_{e}(A_Kx)
	\end{equation*}
	Again, we have
	\begin{equation*}
		V_{e}(x^+) - V_{e}(A_Kx)  = 2(A_Kx)'P e(x) + e(x)'Pe(x) + p'e(x)
	\end{equation*}
	Since $f(x,Kx,0)$ is twice continuously differentiable, there again exist 
	$a_3,a_4>0$ such that
	\begin{equation*}
		V_{e}(x^+) - V_{e}(A_Kx) \leq a_3|x|^3 + a_4|x|^4
	\end{equation*}
	for all $x\in\bbX_f$. Since $\bbX_f$ is bounded, we can choose $\mu\geq 0$ 
	such that
	\begin{equation*}
		a_3|x|^3 + a_4|x|^4 \leq \frac{\mu}{2}c_2|x|^2 \quad \forall \ 
		x\in\bbX_f
	\end{equation*}
	and therefore \cref{eq:termcostdec} holds. 
\end{proof}

\begin{rem}[$\mu=0$]\label{rem:mu}
	Lemma \ref{lem:construct} permits $\mu=0$ if $a_3=a_4=0$, i.e., if 
	$V_e(A_Kx)$ 
	overestimates the value of $V_e(x^+)$ for all $x\in\bbX_f$. As shown in 
	\Cref{s:Example}, we can choose $\mu=0$ for nonlinear systems while still 
	satisfying \cref{as:term}. 
\end{rem}


\section{Example: CSTR}\label{s:Example}
We consider a first-order, irreversible chemical reaction ($A\rightarrow B$) in 
an isothermal CSTR, as 
discussed in \citet{diehl:amrit:rawlings:2010,amrit:rawlings:angeli:2011}. The 
dynamics are 
\begin{align*}
	\frac{dc_A}{dt} & = \frac{q_f}{10}(1+w-c_A) - 0.4c_A \\
	\frac{dc_B}{dt} & = -\frac{q_f}{10}(c_B) + 0.4c_A
\end{align*}
in which $c_A,c_B\in[0,1]$ are the concentration of species $A$, $B$ in the 
reactor, $q_f\in[0,10]$ is the inlet flow rate, and $w\in[-0.2,0.2]$ represents 
a disturbance in the inlet concentration of species $A$. The system 
is discretized with a sample time of $\Delta=0.25$. The economic stage cost is 
\begin{equation*}
	\ell(c_A,c_B,q_f) = -2q_fc_B + (1/2)q_f
\end{equation*}
The optimal steady state for this cost is $c_A=c_B=0.5$ and $q_f=4$. We define 
$x=\begin{bmatrix}
	c_A-0.5 & c_B-0.5
\end{bmatrix}'$ and $u=q_f-4$. Therefore, $\bbX=\{x\mid -0.5\leq x\leq 0.5\}$ 
and $\bbU=\{u\mid -4\leq u \leq 6\}$. 

We also consider the following regularized 
stage cost from \citet{amrit:rawlings:angeli:2011}, which guarantees strict 
dissipativity:
\begin{equation*}
	\ell_d(c_A,c_B,q_f) = -2q_fc_B + (1/2)q_f +0.1(q_f-4)^2 
\end{equation*}
We subsequently compare economic MPC formulations with nondissipative 
$\ell(\cdot)$ and dissipative $\ell_d(\cdot)$ stage costs. 

We construct the terminal cost and constraint according to the approach in 
Section 4 with $K=\begin{bmatrix} -0.012 & -0.037\end{bmatrix}$.\footnote{We 
linearize the continuous time differential equation first and then convert to 
discrete time to give $x^+=Ax+Bu$.} We have that 
$V_s(x)=x'\tilde{P}x$ is a Lyapunov function for the system on $\bbX$, in which 
$\tilde{P}$ satisfies $A_K\tilde{P}A_K - \tilde{P} + I = 0$.  

We define \cref{eq:V_f} with $\mu=0$ (see Remark \ref{rem:mu}). For the purely 
economic 
stage cost 
$\ell(\cdot)$, we have $V_f(x)=x'Px+p'x$ with
\begin{equation*}
	P=\begin{bmatrix}
		-9.47\times 10^{-5} & 4.56\times 10^{-2} \\
		4.56\times 10^{-2} & 4.49\times 10^{-1}
	\end{bmatrix} \quad p = \begin{bmatrix}
	-39.9 \\ -84.1
\end{bmatrix}
\end{equation*} 
For the dissipative stage cost $\ell_d(\cdot)$, we have
\begin{equation*}
	P_d=\begin{bmatrix}
		-5.14\times 10^{-5} & 4.58\times 10^{-2} \\
		4.58\times 10^{-2} & 4.5\times 10^{-1}
	\end{bmatrix} \quad p_d = \begin{bmatrix}
		-39.9 \\ -84.1
	\end{bmatrix}
\end{equation*} 
Note that neither of these terminal cost functions are convex. 
In both cases, we verify that \cref{as:term} holds for 
$\tau=\max_{x\in\bbX}V_s(x)$ and therefore $\bbX_f=\bbX$. 
We choose a horizon of $N=16$ to emphasize the ability of terminal 
costs/constraints to handle problems with short horizons.

\begin{figure}
	\centering
	\includegraphics[width=\linewidth]{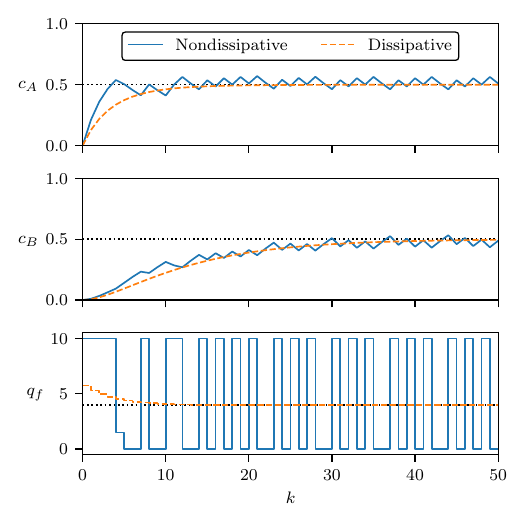}
	\caption{Nominal closed-loop trajectories for CSTR example 
	with nondissipative and dissipative stage costs.}
	\label{fig:econmpc}
\end{figure}

In \cref{fig:econmpc}, we plot the nominal ($w=0$) closed-loop trajectory for 
both stage costs starting from $c_A=c_B=0$. Note that the nondissipative stage 
cost follows a (seemingly) periodic trajectory while the dissipative stage cost 
stabilizes the specified steady state. If frequent actuation of $q_f$ is not 
acceptable, then the dissipative stage cost is preferable. However, if frequent 
actuation of $q_f$ is not a significant issue for the process, then dynamic 
operation via the nondissipative stage cost is economically superior.
 
\begin{figure}
	\centering
	\includegraphics[width=\linewidth]{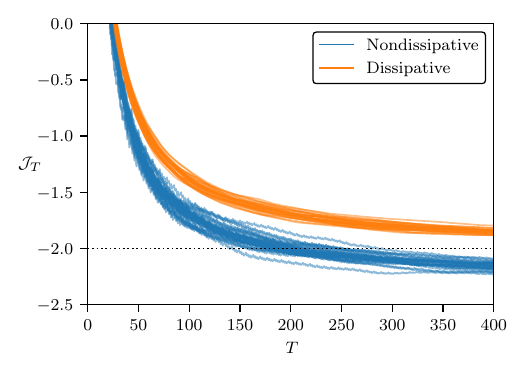}
	\caption{Closed-loop average economic performance for 30 realizations of 
	disturbance trajectory with nondissipative and dissipative 
	stage costs.}
	\label{fig:econmpc_pert}
\end{figure}

In \cref{fig:econmpc_pert}, we plot the closed-loop performance for 30 
realizations of the disturbance trajectory. We define the closed-loop 
performance of these trajectories via the average \textit{economic} 
stage cost:
\begin{equation*}
	\mathcal{J}_T = \frac{1}{T}\sum_{k=0}^{T-1}\ell(x(k),u(k))
\end{equation*}
in which $x(k)=\phi(k;x,\bfw_{\infty})$ and $u(k)=\kappa(x(k))$. Note that we 
evaluate the performance of both economic MPC formulations via only the 
economic (nondissipative) stage cost. For large $T$, we observe that the 
nondissipative stage cost produces better economic performance than the 
dissipative stage cost even when subject to perturbations. 

\section{Future extensions}

We plan to extend these results to suboptimal economic MPC algorithms that can 
be deployed online with limited computation time for high dimensional 
systems. Moreover, we can extend these results to time-varying systems and 
reference trajectories to address a wider class of problems. Many time-varying 
economic MPC applications do not require stability of a target reference 
trajectory and are instead primarily concerned with economic 
performance, e.g., HVAC or production scheduling.

\bibliography{main}



\end{document}